\documentclass[12pt]{article}
\usepackage{amsmath}
\usepackage{amsfonts}
\usepackage{amssymb}
\usepackage{longtable}

\newcommand{\Real}{\mathbb{R}}
\newcommand{\g}{\mathfrak{g}}
\newcommand{\h}{\mathfrak{h}}
\newcommand{\so}{\mathfrak{so}}
\newcommand{\simil}{\mathfrak{sim}}
\newcommand{\f}{\mathfrak{f}}
\newcommand{\R}{\mathcal{R}}
\newcommand{\Ric}{\mathop{{\rm Ric}}\nolimits}
\newcommand{\tRic}{\mathop{\widetilde{\rm Ric}}\nolimits}
\newcommand{\pr}{\mathop\text{\rm pr}\nolimits}
\newcommand{\id}{\mathop\text{\rm id}\nolimits}
\newcommand{\tr}{\mathop\text{\rm tr}\nolimits}

\renewcommand{\d}{\textrm{d}}

\newtheorem{lemma}{Lemma}
\newtheorem{theorem}{Theorem}

\makeindex

\begin{document}

\title{Some applications of the Lorentzian holonomy algebras}
\author{Anton S. Galaev}

\maketitle

\begin{abstract} It is shown how one can apply the classification of the holonomy algebras of Lorentzian manifolds
to solve some problems. In particular, a new proof to the
classification of Lorentzian manifolds with recurrent curvature
tensor
 is given; the classification of two-symmetric Lorentzian manifolds is explained;  conformally recurrent Lorentzian manifolds are classified; recurrent symmetric
 bilinear forms on Lorentzian manifolds are described.
\end{abstract}


\section{Introduction}\label{sec:1}

While the classification of the Riemannian holonomy algebras is a classical result that has many applications
both to geometry and physics, see e.g. \cite{Besse,Joyce07}, the classification of the Lorentzian holonomy algebras
has been achieved  only recently \cite{ESI,Leistner}. We review it in Section \ref{SecClas}.
The holonomy algebra  of a pseudo-Riemannian manifold is an important invariant of the Levi-Civita connection.
It provides information about parallel and recurrent tensors on the manifold. Using that property, we
solve some problems in Lorentzian geometry. As the first illustration, in Section \ref{SecRec}
 we give a new and modern proof to the
classification of Lorentzian manifolds $(M,g)$  with recurrent
curvature tensor $R$, i.e. satisfying \begin{equation}\label{Rrec}
\nabla_X R=\theta(X)R\end{equation} for all vector fields $X$ and
a 1-form $\theta$. Originally this classification is achieved in
\cite{WalkerRec}. In Section \ref{SecSym} we discuss Lorentzian
symmetric spaces. As a new result, in Section \ref{SecWrec} we
obtain a classification of Lorentzian manifolds with recurrent
conformal Weyl tensor $W$. This generalizes result from
\cite{DR77,DR09} that gives classification of Lorentzian manifolds
with parallel $W$. In Section \ref{SecTwoSym} we explain the
result from \cite{A-G} about the classification of two-symmetric
Lorentzian manifolds $(M,g)$, i.e. manifolds satisfying the
condition
\begin{equation}\label{RTwoSym} \nabla^2 R=0,\quad \nabla R\neq
0.\end{equation} In Section \ref{SecRecBF} we study recurrent
symmetric  bilinear forms on Lorentzian manifolds.

\section{Holonomy algebras; parallel and recurrent tensor fields}

Let $(M,g)$ be a connected pseudo-Riemannian manifold of signature
$(r,s)$. {\it The holonomy group} $G_x$ of $(M,g)$ at a point
$x\in M$ is the Lie group that consists of the pseudo-orthogonal
transformations given by the parallel transports along all
piecewise smooth loops at the point $x$. It can be identified with
a Lie subgroup of the pseudo-orthogonal Lie group ${\rm
O}(r,s)={\rm O}(T_xM,g_x)$. The corresponding subalgebra $\g_x$ of
$\so(r,s)=\so(T_xM,g_x)$ is called {\it the holonomy algebra} of
$(M,g)$ at the point $x\in M$.

{\it The Ambrose-Singer Theorems} states that the holonomy algebra
$\g_x$ is spanned by the following endomorphisms of $T_xM$:
$$\tau^{-1}_\gamma\circ R_y(\tau_\gamma X,\tau_\gamma Y)\circ
\tau_\gamma, $$ where  $\gamma$ is a piecewise smooth curve
starting at the point $x$ with an end-point $y \in M$, and $X,Y\in
T_xM$.

Since the manifold $M$ is connected, the holonomy groups (holonomy
algebras) of $(M,g)$ at different points are isomorphic, and one
may speak about the holonomy group $G\subset {\rm O}(r,s)$ (the
holonomy algebra $\g\subset\so(r,s)$) of $(M,g)$.

Recall that a tensor field $T$ on $(M,g)$ is {\it parallel} if
$\nabla T=0$, where $\nabla$ is the Levi-Civita connection; $T$ is
called {\it recurrent} if $\nabla T=\theta\otimes T$ for a 1-form
$\theta$.

 The {\it fundamental principle} \cite{Besse}
states that there exists a one-to-one correspondents between
parallel tensor fields $T$ on $M$ and tensors $T_0$ of the same
type at $x$ preserved by the holonomy group (more precisely, by
its  tensor extension of its representation). Similarly, there
exists a one-to-one correspondents between rank one parallel
subbundles of a tensor bundle  over  $M$ and 1-dimensional
subspaces of the space of tensors  of the same type at $x$
preserved by the holonomy group. Any section of a rank one
parallel subbundle of a tensor bundle is a recurrent tensor field.
Conversely, any non-vanishing recurrent tensor field defines such
parallel subbundle.

If the manifold $M$ is simply connected, then the holonomy group
is connected and it is uniquely defined by the holonomy algebra.
Then parallel and recurrent tensors may be described in terms of
the holonomy algebra.

\section{Classification of the Lorentzian holonomy algebras}\label{SecClas}

Here we review results from \cite{Leistner,ESI}. Let $(M,g)$ be a
simply connected Lorentzian manifold of dimension $n+2$, $n\geq
0$. Fix a point $x\in M$. The tangent space $(T_xM,g_x)$ can be
identified with the Minkowski space $(\Real^{1,n+1},\eta)$. Then
the holonomy algebra  $(M,g)$ at the point $x$ is identified with
a subalgebra $\g\subset\so(1,n+1)$.

We may assume that the holonomy algebra $\g\subset\so(1,n+1)$ of
$(M,g)$ is weakly irreducible, i.e. it does not preserve any
non-degenerate proper vector subspace of the tangent space.
Indeed, if $\g\subset\so(1,n+1)$ is not weakly irreducible, then
by the Wu Theorem,   $(M,g)$ is at least locally can be decomposed
into a product of a Lorentzian manifold and a Riemannian manifold,
see e.g. \cite{ESI}. Thus we  assume that $(M,g)$ is locally
indecomposable.
 If $\g\subset\so(1,n+1)$ is irreducible, then $\g=\so(1,n+1)$.
Suppose that $\g\subset\so(1,n+1)$ is not  irreducible, then $\g$
preserves an isotropic line in $\Real^{1,n+1}$.

The Lie algebra $\so(1,n+1)$ can be identified with the space of
bivectors $\Lambda^2\Real^{1,n+1}$ in such a way that $$(X\wedge
Y)Z= \eta (X,Z)Y- \eta (Y,Z)X.$$ Let $p\in\Real^{1,n+1}$ be an
isotropic vector. Fix an isotropic vector $q\in\Real^{1,n+1}$ such
that $\eta(p,q)=1$. Let $E$ be the orthogonal complement to $\Real
p\oplus\Real q$, then $E\simeq \Real^n$ is an Euclidean space and
we get
$$\Real^{1,n+1}=\Real p\oplus E\oplus\Real q.$$ Denote by
$\simil(n)$ the maximal subalgebra in $\so(1,n+1)$ preserving the
isotropic line $\Real p$, then it holds $$\simil(n)=\Real p\wedge
q+\so(n)+p\wedge E,$$ here $\so(n)=\so(E)\simeq \wedge^2E$. Any
weakly irreducible not irreducible subalgebra
$\g\subset\so(1,n+1)$ preserves an isotropic line in
$\Real^{1,n+1}$, hence $\g$ is conjugated to a subalgebra of
$\simil(n)$. The weakly irreducible Lorentzian holonomy algebras
$\mathfrak{g}\subset\simil(n)$ are
the following: \\
(type I) $\g=\Real p\wedge q+\h+p\wedge E$;\ \ \
(type II) $\g=\h+p\wedge E$;\\
(type III) $\g=\{\varphi(A)p\wedge q+A|A\in\h\}+p\wedge E$;\\
(type IV) $\g=\{A+p\wedge \psi(A) |A\in\h\}+p\wedge E_1$, \\
where  $\h\subset\so(n)$ is a Riemannian holonomy algebra,
$\varphi:\h\to\Real$ is a linear map that is zero on the commutant
$[\h,\h]$; for the last algebra, $E=E_1\oplus E_2$ is an
orthogonal decomposition, $\h$ annihilates $E_2$, i.e. $\h\subset
\so(E_1)$, and $\psi:\h\to E_2$ is a surjective linear map that is
zero on  $[\h,\h]$.

\section{The spaces of curvature tensors}\label{SecCurv}

We will need the following algebraic results. Let $(W,\eta)$ be a
pseudo-Euclidean space and $\mathfrak{f}\subset\so(W)$ be a
subalgebra.  The vector  space
$$\R(\mathfrak{f})=\{R\in \Lambda^2 W^*\otimes\mathfrak{f}|R(X,Y)Z+R(Y,Z)X+R(Z,X)Y=0,\, X,Y,Z\in W\}$$ is called {\it the space of algebraic curvature
tensors of type} $\mathfrak{f}$. The space $\R(\f)$ is an
$\f$-module with the action $$(\xi\cdot R)(X,Y)=[\xi,R(X,Y)]-R(\xi
X,Y)-R(X,\xi Y),\quad\xi\in\f,\, R\in\R(\f).$$  From the
Ambrose-Singer Theorem it follows that  if
$\mathfrak{f}\subset\so(W)$ is the holonomy algebra of a
pseudo-Riemannian manifold $(N,h)$, then the values of the
curvature tensor of $(N,h)$ belong to $\R(\mathfrak{f})$ and
$$\mathfrak{f}={\rm span}\{R(X, Y)|R\in\R(\mathfrak{f}),\,X,Y\in W\},$$ i.e. $\mathfrak{f}$ is spanned
by the images of the elements $R\in\R(\mathfrak{f})$.

The spaces $\R(\g)$ for Lorentzian holonomy algebras
$\g\subset\simil(n)$ are found in \cite{Gal1,onecomp}. Let e.g.
$\g=\Real p\wedge q+\h+p\wedge E$. For the subalgebra
$\h\subset\so(n)$ define the space
$$\mathcal{P}(\h)=\{P\in E^*\otimes\h|\eta(P(X)Y,Z)+\eta(P(Y)Z,X)+\eta(P(Z)X,Y)=0,\,
X,Y,Z\in E\}.$$ Any $R\in\R(\g)$ is uniquely given by
$$\lambda\in\mathbb{R},\ \vec v\in E,\
P\in\mathcal{P}(\mathfrak{h}),\
R_0\in\mathcal{R}(\mathfrak{h}),\text { and }T\in {\rm
End}(E)\text{ with }T^*=T$$ in the following way:
\begin{align*}
R(p,q)=& -\lambda p\wedge q-p\wedge \vec v,\,\,\, R(X,Y)=R_0(X,Y)-p\wedge (P(Y)X-P(X)Y),\\
R(X,q)=&-g(\vec{v},X)p\wedge q+P(X)-p\wedge T(X),\qquad R(p,X)=0
\end{align*}  for all $X,Y\in E$.
For the algebras $\g$ of the other types, any $R\in\R(\g)$ can be
given in the same way and by the condition that $R$ takes values
in $\g$. For example, $R\in\R(\h+p\wedge E)$ if and only if
$\lambda=0$ and $\vec v=0$.

\section{Walker metrics and pp-waves}
Consider the local form of a Lorentzian manifold $(M,g)$ with the
holonomy algebra $\g\subset\simil(n)$. Since $\g$  preserves an
isotropic line of the tangent space,  $(M,g)$ locally admits a
parallel distribution of isotropic lines. Locally there exist the
so called Walker coordinates $v,x^1,...,x^n,u$ and the metric $g$
has the form
\begin{equation}\label{Walker} g=2\d v\d u+h+2A\d u+H (\d
u)^2,\end{equation} where $h=h_{ij}(x^1,...,x^n,u)\d x^i\d x^j$ is
an $u$-dependent family of Riemannian metrics,  $A=A_i(x^1, \ldots
, x^n,u)\d x^i$ is an $u$-dependent family of one-forms, and $H$
is a local function on $M$, see e.g. \cite{ESI}. The vector field
$\partial_v$ defines the parallel distribution of isotropic lines.

Consider the fields of frames $$p=\partial_v,\quad
X_i=\partial_i-A_i\partial_v,\quad
q=\partial_u-\frac{1}{2}H\partial_v$$ and the distribution $E={\rm
span}\{X_1,...,X_n\}.$ At each point $m$ of the coordinate
neighborhood we get the decomposition
$$T_mM=\Real p_m\oplus E_m\oplus\Real q_m,$$ hence the value $R_m$
of the curvature tensor can be expressed in terms of some
$\lambda_m$, $\vec v_m$, $R_{0m}$, $R_m$ and $T_m$ as above. The
space $E_m$ is isomorphic to the tangent space of a Riemannian
manifolds with a metric from the family $h$, then $R_{0}$ is
defined by the curvature tensor of the family of the Riemannian
metrics~$h$.

It is known \cite{ESI} that the holonomy algebra of the manifold
$(M,g)$ is contained in $p\wedge E\subset\simil(n)$ if and only if
the metric can be locally written in the form
 \begin{equation}\label{pp-wave} g=2\d v\d u+\sum_{i=1}^n(\d x^i)^2+H (\d u)^2,\quad\partial_v H=0.\end{equation}
 Such spaces are called pp-waves.

\section{Lorentzian manifolds with recurrent curvature tensor}\label{SecRec}

In this section we consider Lorentzian manifolds $(M,g)$  with recurrent curvature tensor $R$, i.e.
satisfying  \eqref{Rrec}. Note that for Riemannian
manifolds \eqref{Rrec}, implies $\theta=0$, i.e. the manifold is
locally symmetric \cite{Kai}.

 Many facts about recurrent spaces, or more generally
about $r$-recurrent spaces, and a long list of literature on this
topic can be found in the fundamental review of Kaigorodov
\cite{Kai}. There is a recent review by Senovilla \cite{Sen08},
where similar problems are considered.

In this section we give a new proof to the following theorem
proven by Walker in~\cite{WalkerRec}.

\begin{theorem}\label{ThRrec}
Let $(M,g)$ be a Lorentzian manifold of dimension $n+2\geq 3$.
Then $(M,g)$ is recurrent and not locally symmetric if and only if
in a neighborhood of each point of $M$ there exist coordinates
$v,x^1,...,x^n,u$ such that one of the following holds:
\begin{itemize}\item[I.]  there exists a function
$H(x^1,u)$ such that \begin{equation}\label{gRrec0} g=2\d v\d
u+\sum_{i=1}^n(\d x^i)^2+H(x^1,u)(\d u)^2.\end{equation}

\item[II.] There exist real numbers $\lambda_1,...,\lambda_n$ with $|\lambda_1|\geq \cdots\geq
|\lambda_n|$, $\lambda_2\neq 0$, and a function \\
$F:U\subset\Real\to\Real$   such that
\begin{equation}\label{gRrec}
g=2\d v\d u+\sum_{i=1}^n(\d x^i)^2+F(u)\lambda_i(x^i)^2(\d
u)^2.\end{equation}

Moreover, for some system of coordinates $\partial^2_{1}H$ is not
constant or $\frac{\d F}{\d u}\neq 0$.
\end{itemize}
\end{theorem}

The form of the metric may change from one system coordinates to
another, i.g. it can be flat for some systems of coordinates.
Examples of such spaces can be constructed taking  the metrics of
the form \eqref{gRrec} with $F(u)=0$ if $|u|\geq\epsilon$ for some
$\epsilon>0$, any such metric is flat on the spaces
$\{(v,x^1,\dots x^n,u)|\,|u|\geq\epsilon\}$, hence me may  glue
these metrics on such flat spaces. In this example the function
$F(u)$ is not analytic. Theorem \ref{ThRrecAn} below states that
if the manifold $(M,g)$ is analytic, then the metric is the same
for all systems of coordinates.

Note that the local metric \eqref{gRrec} is symmetric if and only
 $F$ is a constant, i.e. $\frac{\d F(u)}{\d u}=0$. In this case
we get the so called Cahen-Wallach space \cite{C-W}. Next, the
local metric \eqref{gRrec} is two-symmetric, i.e. $\nabla^2R=0$,
if and only if $\frac{\d^2 F(u)}{(\d u)^2}=0$, see Section
\ref{SecTwoSym} below. Finally, it is conformally flat if and only
if $\lambda_1=\cdots=\lambda_n$ \cite{GalConf}.

\subsection{Proof of Theorem \ref{ThRrec}}

First we reduce the problem to the case when $(M,g)$ is locally
indecomposable.

\begin{lemma} Let $(M,g)$ be a recurrent and not locally symmetric Lorentzian manifold.
Suppose that $(M,g)$ is locally decomposable, i.e. each point of
$M$ has an open neighborhood $U$ such that $(U,g|_{U})$ is
isometric to the product of a Lorentzian manifold $(M_1,g_1)$ and
a Riemannian manifold $(M_2,g_2)$. If $\nabla R|_{U}\neq 0$, then
$(M_1,g_1)$ is recurrent and $(M_2,g_2)$ is flat. If $\nabla
R|_{U}=0$, then both  $(M_1,g_1)$  and $(M_2,g_2)$ are locally
symmetric. \end{lemma}

{\bf Proof.}  Since $(U,g|_{U})=(M_1\times M_2,g_1+g_2)$, for the
corresponding curvature tensors and their covariant derivatives it
holds
$$R|_{U}=R_1+R_2,\quad \nabla R|_{U}=\nabla R_1+\nabla R_2.$$
Suppose that $\nabla R|_{U}\neq 0$. Restricting the equality
$\nabla R=\theta\otimes R$ to $(M_2,g_2)$, we get $\nabla
R_2=\theta|_{M_2}\otimes R_2$. Since $(M_2,g_2)$ is a Riemannian
manifold, $\theta|_{M_2}=0$. Let $X_1\in\Gamma(TM_1)$ and
$X_2,Y_2\in\Gamma(TM_2)$, then
\begin{multline*}0=\nabla_{X_1}R_1(X_2,Y_2)+\nabla_{X_1}R_2(X_2,Y_2)\\=\theta(X_1)R_1(X_2,Y_2)+\theta(X_1)R_2(X_2,Y_2)=\theta(X_1)R_2(X_2,Y_2).\end{multline*}
Since $\theta |_{U}\neq 0$, $R_2=0$. This proves the lemma. $\Box$

The condition \eqref{Rrec} implies that for any point $m\in M$,
the holonomy algebra $\g_m$ of $(M,g)$ preserves the line $\Real
R_m\subset\R(\g_m)$ in the space of possible values of the
curvature tensor at the point $m$.

The only possible irreducible holonomy algebra of $(M,g)$ is the
Lorentzian Lie algebra $\so(1,n+1)$ \cite{ESI}. Form the results
of \cite{Al} it follows that the only line preserved by
$\so(1,n+1)$ in the space $\R(\so(1,n+1))$ consists of curvature
tensors defined by the scalar curvature. Consequently the manifold
is Einstein and  locally symmetric. Hence the holonomy algebra of
$(M,g)$ is weakly irreducible and not irreducible and it is
contained in $\simil(n)$.

The condition that the holonomy algebra $\g_m$ at the point $m\in
M$ preserves the line $\Real R_m\subset\R(\g_m)$ can be expressed
as
$$\xi\cdot R_m=\mu(\xi)R_m,\quad\xi\in\g_m,$$ where $\mu:\g_m\to\Real$ is a linear
map. Let e.g. $\g_m=\Real p_m\wedge q_m+\h+p_m\wedge E_m$. As the
$\h$-module, the space $\R(\g_m)$ admits the decomposition
$$\R(\g)=\Real\oplus E_m\oplus\R(\h)\oplus\mathcal{P}(\h)\oplus\odot^2E_m.$$
The space $\mathcal{P}(\h)$ does not contain any $\h$-invariant
one-dimensional subspace \cite{onecomp}, hence $P_m=0$. For
$X,Y,Z\in E_m$ it holds \begin{multline*}\mu(p_m\wedge
Z)R_{0m}(X,Y)=\mu(p_m\wedge Z)R_{m}(X,Y)=((p_m\wedge Z)\cdot
R_{m})(X,Y)\\=[p_m\wedge Z,R_{0m}(X,Y)]=-p_m\wedge
R_{0m}(X,Y)Z.\end{multline*} This implies $R_{0m}=0$. Thus over
the current coordinate neighborhood it holds $R_0=0$ and $P=0$.
The same can be shown for the other possible holonomy algebras. We
get $R(p^\bot,p^\bot)=0$. In \cite{Le06} it is proved that in this
case the coordinates can be chosen in such a way that
\begin{equation} g=2\d v\d u+\sum_{i=1}^n(\d x^i)^2+H (\d u)^2.\end{equation} In
particular, $\h=0$ and either $\g_m=p_m\wedge E_m$, or $\g_m=\Real
p_m\wedge q_m+p_m\wedge E_m$. Consider these two cases.

{\bf Case 1.} Suppose that $\g_m=p_m\wedge E_m$. Then $\partial_v
H=0$. In \cite{A-G} it is shown that the covariant curvature
tensor and its covariant derivative have the form \begin{align}
\label{Rpp-wave}  R&= \frac{1}{2}(\partial_i\partial_jH)(q' \wedge
e^i \vee q' \wedge e^j ),\\\label{nablaRpp-wave} \nabla
R&=\frac{1}{2}(\partial_k\partial_i\partial_jH)e^k\otimes(q'
\wedge e^i \vee q' \wedge e^j
)+\frac{1}{2}(\partial_u\partial_i\partial_jH)q'\otimes(q' \wedge
e^i \vee q' \wedge e^j ),\end{align} where $e^i=\d x^i$ and $q'=\d
u$. The condition \eqref{Rrec} is equivalent to
$$\partial_k\partial_i\partial_jH=\theta_k\partial_i\partial_jH,\quad
\partial_u\partial_i\partial_jH=\theta_u\partial_i\partial_jH,$$ where $\theta_k=\theta(\partial_k)$
and $\theta_u=\theta(\partial_u)$. If $\partial_i\partial_j H\neq
0$ for some $i,j$ on some open subspace, then
$$\theta_k=\partial_k\ln|\partial_i\partial_j H|,\quad \theta_u=\partial_u\ln|\partial_i\partial_j
H|,$$ i.e. $\d\theta=0$ and there exists a function $f$ such that
$\theta=\d f$. We get
$$\partial_k(\ln|\partial_i\partial_j H|-f)=\partial_u(\ln|\partial_i\partial_j H|-f)=0,$$
 i.e. $$\ln|\partial_i\partial_j H|=f+c_{ij},\quad c_{ij}\in \Real,\quad
 c_{ij}=c_{ji}.$$ Thus, $$\partial_i\partial_jH=e^fC_{ij},\quad C_{ij}=e^{c_{ij}}.$$
Consider the new coordinates $$\tilde v=v,\quad\tilde
x^i=a^i_jx^j,\quad \tilde u=u,$$ where $a^i_j$ is an orthogonal
matrix. With respect to these coordinates the metric $g$ takes the
same form and it holds $$\tilde \partial_i\tilde \partial_j\tilde
H=e^{\tilde f} a^r_ia^s_jC_{rs}.$$ The orthogonal transformation
$a^j_i$ can be chosen in such a way that the matrix $\tilde
C_{ij}=a^r_ia^s_jC_{rs}$ is diagonal with the diagonal elements
$\lambda_1,...,\lambda_n$. Assume that
$|\lambda_1|\geq\cdots\geq|\lambda_n|$. Thus  it holds
$$\partial_i\partial_jH=e^f\delta_{ij}\lambda_i,\quad \lambda_i\in\Real.$$

If $\lambda_2=\cdots=\lambda_n=0$, then
$$H=F(x^1,u)+\sum_{i=2}^nG_i(u)x^i.$$
Consider the new coordinates given by the inverse transformation
\begin{equation}\label{trans1} u =\tilde u, \quad x^i = \tilde x^i+b^i(\tilde
u), \quad v = \tilde v  -\sum_j \frac{\d b^j(\tilde u)}{\d \tilde
u}\tilde x^j\end{equation} such that $2\frac{\d^2b^j(u)}{(\d
u)^2}=G_j(u)$ and $b^1(u)=0$. With respect to the new coordinates
it holds $H=F(x^1,u)$ and we obtain the Case I of the formulation
of the theorem.

Suppose that $\lambda_2\neq 0$. From the above we get that if
$i\neq j$, then $\partial_i\partial_jH=0$, i.e. $H$ is of the form
$H=\sum_i H_i(x^i)$, and $\frac{d^2 H_i}{(d x^i)^2}=e^f\lambda_i$.
Taking $i=1,2$ and differentiating the last  equality with respect
to $\partial_j$, we get $\partial_jf=0$, i.e. $f$ depends only on
$u$. Now it is clear that
$$H=\frac{1}{2}e^{f(u)}\lambda_i(x^i)^2+B_i(u)x^i+K(u).$$
Let $F(u)=\frac{1}{2}e^{f(u)}$. From the results of \cite{A-G} it
follows that the coordinates can be chosen in such a way that
$H=F(u)\lambda_i(x^i)^2.$

{\bf Case 2.} Suppose that $\g_m=\Real p_m\wedge q_m+p_m\wedge
E_m$. The curvature tensor $R_m$ is given by the elements
$\lambda_m$, $\vec v_m$ and $T_m$. It holds
\begin{multline*}\mu(p_m\wedge q_m)(-\lambda_m p_m\wedge
q_m-p_m\wedge \vec v_m)\\=((p_m\wedge q_m)\cdot
R_m)(p_m,q_m)=[p_m\wedge q_m,R(p_m,q_m)]=p_m\wedge \vec
v_m,\end{multline*} hence $$\mu(p_m\wedge q_m)\lambda_m=0,\quad
(\mu(p_m\wedge q_m)+1)\vec v_m=0.$$ Similarly, $((p_m\wedge
q_m)\cdot R_m)(X,q_m)=\mu(p_m\wedge q_m)R_m(X,q_m)$, for $X\in
E_m$, implies
$$\mu(p_m\wedge q_m)\vec v_m=0,\quad (\mu(p_m\wedge
q_m)+1)T_m=0.$$ In the same way, using an element $p_m\wedge
X\in\g_m$, we get \begin{multline*}\mu(p_m\wedge
X)\lambda_m=0,\quad \mu(p_m\wedge X)\vec v_m=\lambda_m X,\\
\mu(p_m\wedge X)\vec v_m=0,\quad g(\vec v_m,Y)X=\mu(p_m\wedge
X)T_m(Y).\end{multline*} The obtained equalities imply $\vec
v_m=0$ and $\lambda_m=0$. Consequently, over this coordinate
neighborhood, $\lambda=0$ and $\vec v=0$. This shows that this
coordinate neighborhood is the same as in Case 1.
 $\Box$

\subsection{The case of analytic $(M,g)$}

Suppose that $(M,g)$ is analytic. In this case, Theorem
\ref{ThRrec} can be reformulated in the following way:

\begin{theorem}\label{ThRrecAn}
Let $(M,g)$ be an analytic  Lorentzian manifold of dimension
$n+2\geq 3$. Then $(M,g)$ is recurrent and not locally symmetric
if and only if one of the following holds:
\begin{itemize}\item[I.] In a neighborhood of each point of $M$
there exist coordinates $v,x^1,...,x^n,u$ and a function
$H(x^1,u)$ such that \begin{equation}\label{gRrec0An} g=2\d v\d
u+\sum_{i=1}^n(\d x^i)^2+H(x^1,u)(\d u)^2,\end{equation} and
$\partial^2_{1}H$ is not constant for some system of coordinates.

In this case if $n\geq 2$, then the manifold is locally a product
of the three-dimensional recurrent Lorentzian manifold with the
coordinates $v,x^1,u$ and of the flat Riemannian manifold with the
coordinates $x^2,...,x^n$.

\item[II.] There exist real numbers $\lambda_1,...,\lambda_n$ with $|\lambda_1|\geq \cdots\geq
|\lambda_n|$, $\lambda_2\neq 0$, and an analytic function
$F:U\subset\Real\to\Real$ with $\frac{\d F}{\d u}\neq 0$, and in a
neighborhood of each point of $M$ there exist coordinates
$v,x^1,...,x^n,u$  such that \begin{equation}\label{gRrecAn} g=2\d
v\d u+\sum_{i=1}^n(\d x^i)^2+F(u)\lambda_i(x^i)^2(\d
u)^2.\end{equation}

The manifold $(M,g)$ is locally indecomposable if and only if all
$\lambda_i$ are non-zero.

If for some $r$ ($2\leq r<n$) it holds $\lambda_r\neq 0$ and
$\lambda_{r+1}=\cdots=\lambda_n=0$, then  $(M,g)$ is locally  a
product of the recurrent Lorentzian manifold with the coordinates
$v,x^1,\dots,x^r,u$ and of the flat Riemannian manifold with the
coordinates $x^{r+1},...,x^n$.
\end{itemize}
\end{theorem}

In particular, the theorem states that in the second case the
metric is the same in each coordinate neighborhood.

{\bf Proof.}  Suppose that a point $m$ belongs to two coordinate
neighborhoods with the coordinates $v,x^1,\dots,x^n,u$ and $\tilde
v,\tilde x^1,\dots,\tilde x^n,\tilde u$. Suppose that for the
first system of coordinates it holds $H=F(u)\lambda_i(x^i)^2$,
$\lambda_1,\lambda_2\neq 0$, and $\frac{\d F}{\d u}\neq 0$, i.e.
the metric restricted to the first coordinate neighborhood is not
flat. If in the second coordinate system the metric is flat, then
on the intersection of the coordinate domains it holds $\frac{\d
F}{\d u}= 0$. Since $F$ is analytic, this implies $\frac{\d F}{\d
u}= 0$ for all points of the first coordinate neighborhood and we
get a contradiction (this is the only place, where we use the
analyticity). Since the metric restricted to the second coordinate
neighborhood is not flat, the parallel vector field $\tilde
\partial_v$ is defined up to a constant and we may assume that
$\tilde \partial_v=\partial_v$. Then the transformation of
coordinates must have the form $$ u =\tilde u+c, \quad x^i =
a^i_j\tilde x^j+b^i(\tilde u), \quad v = \tilde v -\sum_j
a^j_i\frac{\d b^j(\tilde u)}{\d \tilde u}\tilde x^i+d(\tilde u).$$
where $c\in\Real$, $a^j_i$ is an orthogonal matrix, and
$b^i(\tilde u)$, $d(\tilde u)$ are some functions of
 $\tilde u$ \cite{A-G}. Clearly, the metric written in the second coordinate system
 can not be as in Case I of Theorem \ref{ThRrec}, i.e. it holds
 $$\tilde H=\tilde F(\tilde u)\tilde\lambda_i(\tilde x^i)^2.$$
 Note that $$\tilde F(\tilde u) \delta_{ij}\tilde\lambda_i=F(\tilde
 u+c)\delta_{kl}\lambda_ka^k_ia^l_j.$$
Since the matrix $a_i^j$ is orthogonal, after some change
$$(F(\tilde u),\tilde \lambda_i)\mapsto \left(\frac{1}{C}F(\tilde u),C\tilde
\lambda_i\right),\quad C\neq 0,$$ we obtain
$\tilde\lambda_i=\lambda_i$ and $\tilde F(\tilde u)=F(\tilde
u+c)$. After the transformation $\tilde u\mapsto \tilde u-c$, we
get $\tilde F=F$. This proves the theorem. $\Box$

\section{Lorentzian symmetric spaces}\label{SecSym}

Classification of simply connected Riemannian symmetric spaces is
a classical result of \'Elie Cartan \cite{Besse}. Simply connected
Lorentzian symmetric spaces are classified by Cahen and Wallach
\cite{C-W,Cahen98}. Here we show how the last result can be
reproved using the holonomy theory. It is well-known that a simply
connected pseudo-Riemannian symmetric space is uniquely defined by
the pair $(\g,R)$, where $\g$ is its holonomy algebra and $R$ is
its curvature tensor at  a fix point. Such pair satisfies
$R\in\R(\g)$, $\g$ annihilates $R$ and the image of $R$ coincides
with  $\g$. Note that $R$ can be defined up to a positive
constant.
 Now we describe such pairs  for $\g\subset\so(1,n+1)$.
It is enough to consider indecomposable spaces, i.e. we may assume
that    $\g\subset\so(1,n+1)$ is weakly irreducible. If
$\g=\so(1,n+1)$, then such $R$   constitute a one-dimensional
space without the zero. The connected components of this space
define de Sitter and  Anti de Sitter spaces.  These are the only
indecomposable simply connected Lorentzian symmetric spaces with
semisimple isometry group (equivalently, with reductive holonomy
algebra). Now we suppose that $\g\subset\simil(n)$. In the same
way as we did in Section \ref{SecRec}, we conclude that
$\g=p\wedge E$ and $R$ is given by a symmetric endomorphism $T$ of
$E$ (in notation of Section \ref{SecCurv}).  Thus such pair
corresponds to a pp-wave \eqref{pp-wave}. Equation
\eqref{nablaRpp-wave} shows that $\partial_k\partial_i\partial_j
H=\partial_u\partial_i\partial_j H=0$, i.e.
$H=a_{ij}x^ix^j+b_i(u)x^i+c(u)$. Changing the coordinates, we get
$H=\sum_i\lambda_i (x^i)^2$, $\lambda_i\in\Real$. Metric
\eqref{pp-wave} with such $H$ is defined on $\Real^{n+2}$ and it
is complete. These symmetric spaces are called the Cahen-Wallach
spaces.

\section{The Weyl conformal curvature tensor of a Walker metric}

Below we will need the expression for the Weyl tensor $W$ of a
Walker metric in terms of notations of Section \ref{SecCurv}. This
expression is obtained in \cite{GalConf}. It holds, $$W=R+R_L,$$
where $R_L$ is defined as follows:
\begin{align}
R_L(p,X)&=\frac{1}{n}p\wedge\left(\Ric(h)+\frac{(n-1)\lambda-s_0}{n+1}\id\right)
X,\\
R_L(p,q)&=\frac{1}{n}\left(\frac{2n\lambda-s_0}{n+1}p\wedge
q+p\wedge (\vec{v}-\tRic P)\right),\\
R_L(X,Y)&=\frac{1}{n}\left(p\wedge((X\wedge Y)(\vec v-\tRic P))\right.\\
\nonumber &\left.+\left(\Ric(h)-\frac{s}{2(n+1)}\right)X\wedge
Y+X\wedge
\left(\Ric(h)-\frac{s}{2(n+1)}\right)Y\right),\\
R_L(X,q)&=\frac{1}{n}\Big((\tr T)p\wedge X+ g(X,\vec{v}-\tRic
P)p\wedge q+X\wedge(\vec{v}-\tRic
P)\\\nonumber&+\left(\Ric(h)+\frac{(n-1)\lambda-s_0}{n+1}\id\right)X\wedge
q\Big),
\end{align}
where $\Ric(h)$ is the Ricci operator of the metric $h$,
$s=2\lambda+s_0$ is the scalar curvature of $g$ and $S_0$ is the
scalar curvature of $h$.
 This expression is used in \cite{GalConf} to find all conformally flat Walker metrics.

\section{Lorentzian manifolds with recurrent and parallel Weyl
tensor}\label{SecWrec}

Conformally symmetric Lorentzian manifolds, i.e. Lorentzian
manifolds with  parallel Weyl tensor $W$ are classified by
Derdzinski and Roter \cite{DR77,DR09}. These spaces are exhausted
by  conformally flat spaces, i.e.  with $W=0$, by locally
symmetric spaces, i.e. with $\nabla R=0$, and by some special
pp-waves. As a generalization of this condition one consider
conformally recurrent spaces, i.e. with recurrent Weyl tensor,
$\nabla W=\theta\otimes W$, see e.g. \cite{Ols,SuhK}. We prove the
following theorem.

\begin{theorem}\label{ThWrec} Let $(M,g)$ be a locally indecomposable Lorentzian manifold of dimension $n+2\geq 4$ with
 a recurrent Weyl tensor $W$, then either $W=0$, or $\nabla R=0$, or
 locally $g$ has the form
$$g=2\d v2\d u+\sum_{i=1}^n (\d x^i)^2+\left(a(u)\sum_{i=1}^n(x^i)^2+F(u)\sum_{i=1}^n\lambda_i(x^i)^2\right)(\d u)^2,$$
where $a(u)$, $F(u)$ are functions, and $\lambda_i\in\Real$,
$\sum_{i=1}^n\lambda_i=0$.
\end{theorem}

Note that the for the above metric it holds $\nabla W=0$ if and
only if $\frac{\d F(u)}{\d u}=0$. In particular, we recover the
result by  Derdzinski and Roter. Next, $\nabla R=0$ if and only if
$\frac{\d a(u)}{\d u}=\frac{\d F(u)}{\d u}=0$. Also, $\nabla
R=\theta\otimes R$ if and only if $a(u)=F(u)$, or $a(u)=0$, or all
$\lambda_i=0$. Finally, $W=0$ if and only if all $\lambda_i=0$ or
$F(u)=0$ \cite{GalConf}.

{\bf Proof of Theorem \ref{ThWrec}.} The proof is very similar to
the proof of Theorem \ref{ThRrec} and we omit some obvious
computations. Suppose that $W\neq 0$ and $\nabla R\neq 0$. Let
$\g\subset\so(1,n+1)$ be the holonomy algebra of $(M,g)$ at a
point $m\in M$. Then $\g$ preserves the line in the space $\R(\g)$
spanned by $W_m$. For $\g=\so(1,n+1)$ this would imply $W_m=0$,
which follows from \cite{Al}. Hence, $\g\subset\simil(n)$. Suppose
that $W_m\neq 0$.

\begin{lemma} The manifold $(M,g)$ is a pp-wave, i.e. $\g=p_m\wedge E_m$. \end{lemma}
{\bf Proof.}  For each $\xi\in\g$ it holds $\xi\cdot
W_m=\mu(\xi)W_m$, where $\mu:\g\to\Real$ is a linear map.

First suppose that $\h=0$. Then either $\g=\Real p_m\wedge
q_m+p_m\wedge E_m$,  or $\g=p_m\wedge E_m$.
   Suppose that $\g=\Real p_m\wedge q_m+p_m\wedge E_m$. Then we may assume that $\lambda_m\neq 0$ or $\vec v_m\neq 0$.  Let $Z\in E_m$ be not proportional to $\vec v_m$. Considering $((p_m\wedge Z)\cdot W_m)(p_m,q_m)$, we obtain $\lambda_m=0$. Considering $((p_m\wedge Z)\cdot W_m)(Z_m,q_m)$, we get $\vec v_m=0$.
   We conclude that $\g=p_m\wedge E_m$.

   Suppose now that $\h\neq 0$. Let $A\in\h$ and let $\xi$ be either $A$, or $A+\varphi(A)p_m\wedge q_m$, or $A+p_m\wedge \psi(A)$ depending on the type of $\g$. Note that any 1-dimensional representation of $\h$ is trivial, consequently, $\mu(\xi)=0$. Using this and considering the projection of $(\xi\cdot W_m)(X,Y)$ to $\so(n)$, where $X,Y\in E_m$, we get $$A\cdot\left(\frac{1}{n}(\cdot\wedge (\vec v_m-\tRic P_m))+P_m\right)=0,$$
where we consider the representation of $\h$ in the space
$\mathcal{P}(\so(n))$. The module $\mathcal{P}(\h)$ never contains
non-zero elements annihilated by $\h$ \cite{onecomp}. If $P_m=0$,
then from the above equality it follows that $\vec v_m=0$.
Otherwise, since $\h\neq 0$, there exists $A\in\h$ such that
$A\cdot P_m\neq 0$. This implies that $0\neq \cdot\wedge A(\vec
v_m-\tRic P_m)\in\mathcal{P}(\h)$. Consequently, $\h=\so(n)$. We
conclude that
$$P_m=\frac{1}{n}(\cdot\wedge (\vec v_m-\tRic P_m)).$$
Applying $\tRic$, we get $\tRic P_m=(1-n)\vec v_m$, and we
conclude that
$$P_m=-\cdot\wedge \vec v_m.$$
This shows that the expression of $W_m$ does not include $\vec
v_m$ and $P_m$. Considering $((p_m\wedge Z)\cdot W_m)(p_m,q_m)$,
$Z\in E_m$, we get that $\Ric(h)_m=c_1\id_{E_m}$, where $c_1$ can
be expressed in terms of $\lambda_m$ and $s_{0m}$. Taking the
trace, we get a relation between  $\lambda_m$ and $s_{0m}$. Taking
$\xi$ as above, using  the equality $\pr_{\so(n)}((\xi\cdot
W)(X,Y))=0$, $X,Y\in E_m$, and expressing $R(h)_m$ in terms of the
Weyl tensor of $h$, $\Ric(h)_m$ and $s_{0m}$, we get another
relation between    $\lambda_m$ and $s_{0m}$. Then we conclude
that  $\lambda_m=s_{0m}=0$, $\Ric(h)_m=0$, and $R(h)_m=0$. Now
$W_m$ depends only on $T_m$. Since $\h\neq 0$ and $R(h)=0$, we may
assume that $P_m\neq 0$, and we have just seen that this implies
$\h=\so(n)$. Taking $A\in\so(n)$ and considering $(A\cdot
W_m)(X,q_m)$, we get $T_mA=AT_m$. The Schur Lemma implies that
$T_m$ is proportional to $\id_{E_m}$ . Thus, $W_m=0$. And we get a
contradiction.  Thus, $\h=0$ and $\g=p_m\wedge E_m$, this proves
the lemma.  $\Box$

Now we should find all functions $H$ such the Weyl tensor of
metric \eqref{pp-wave} of a pp-wave is recurrent, i.e. $\nabla
W=\theta\otimes W$ for a 1-form $\theta$. For $W$ and $\nabla W$
we get the formulas \eqref{Rpp-wave} and \eqref{nablaRpp-wave}
with $\partial_i\partial_j H$ replaced by
$\partial_i\partial_jH-\frac{1}{n}\delta_{ij}\Delta H$, where
$\Delta=\sum_{k=1}^n\partial_k^2$. We obtain the equations
$$\partial_k\left(\partial_i\partial_jH-\frac{1}{n}\delta_{ij}\Delta H\right)=\theta_k\left(\partial_i\partial_jH-\frac{1}{n}\delta_{ij}\Delta
H\right),$$
$$\partial_u\left(\partial_i\partial_jH-\frac{1}{n}\delta_{ij}\Delta
H\right)=\theta_u\left(\partial_i\partial_jH-\frac{1}{n}\delta_{ij}\Delta
H\right),$$ where $\theta_k=\theta(\partial_k)$ and
$\theta_u=\theta(\partial_u)$.  As in Section \ref{SecRec}, we get
$$\partial_i\partial_jH-\frac{1}{n}\delta_{ij}\Delta H=e^fC_{ij}, \quad C_{ij}=C_{ji}\in \Real$$
for all $i,j$. Since  $\Delta$ is invariant with respect to an
orthogonal transformation of the coordinates $x^1,...,x^n$, we may
apply a transformation as in Section \ref{SecRec}, ad we may
assume that $$\partial_i\partial_jH-\frac{1}{n}\delta_{ij}\Delta
H=e^f\delta_{ij}\lambda_i,\quad \lambda_i\in\Real.$$ If $i\neq j$,
then $\partial_i\partial_jH=0$, i.e. $H$ is of the form
$H=\sum_{i=1}^n H_i(x^i,u)$. We get the system of equations
$$\partial^2_iH_i-\frac{1}{n}\sum_{k=1}^n\partial_k^2H_k=e^f\lambda_i.$$
We may view this system as a system of linear equations with
respect to the unknowns $\partial^2_iH_i$. The rank of this system
is equal to $n-1$. Summarizing the equations, we see that if a
solution exists, then  $\sum_{i=1}^n\lambda_i=0$. In this case the
dimension of solutions equals 1, and we have
$\partial_i^2H_i=\frac{1}{2}a+e^f\lambda_i$ for a function $a$.
This implies that both $a$ and $f$ are functions depending only on
$u$. We obtain $$H=a(u)\sum_{i=1}^n(x^i)^2+F(u)\lambda_i(x^i)^2,$$
where $F(u)=\frac{1}{2}e^{f(u)}$ and we assume that the terms
linear in $x^i$ are zero, since we can get read of them using a
transformation of coordinates. The theorem is proved. $\Box$

Results about four-dimensional conformally recurrent Lorentzian
spaces are collected in~\cite[Ch. 35]{St}.

\section{Two-symmetric Lorentzian manifolds}\label{SecTwoSym}

In this section we consider two-symmetric Lorentzian manifolds,
i.e. manifolds satisfying \eqref{RTwoSym}. The following theorem
is proved in \cite{A-G}.

\begin{theorem} Let $(M,g)$ be a locally indecomposable Lorentzian manifold of dimension
$n+2$. Then $(M,g)$ is two-symmetric if and only if locally there
exist coordinates $v,x^1,...,x^n,u$ such that
$$g=2\d v\d u +\sum_{i=1}^n(\d x^i)^2+(H_{ij} u+F_{ij})x^ix^j(\d u)^2,$$
where   $H_{ij}$ is a nonzero diagonal real matrix with the
diagonal elements $\lambda_1\leq\cdots\leq\lambda_n$, and $F_{ij}$
is a  symmetric real matrix.
\end{theorem}

 Detailed investigation
of two-symmetric Lorentzian spaces initiated  Senovilla in
\cite{Sen08}, where it is proven  that any two-symmetric
Lorentzian space admits a parallel isotropic vector field, i.e.
locally the metric has the form \eqref{Walker} with
$\partial_vH=0$.

Now we explain the proof of the above theorem from \cite{A-G}. The
assumption that a Lorentzian manifold $(M,g)$ is two-symmetric
implies that the holonomy algebra $\g\subset\so(1,n+1)$ of $(M,g)$
at a point $m\in M$ annihilates the tensor $\nabla R_m\neq 0$. The
tensor $\nabla R_m$ belongs to the $\g$-module $\nabla\R(\g)$ that
consists of linear maps from $\Real^{1,n+1}$ to $\R(\g)$
satisfying the second Bianchi identity.  The results from
\cite{Str88} show that the space $\nabla \R(\so(1,n+1))$ does not
contain any non-zero element annihilated by $\so(1,n+1)$. Hence,
$\g$ can not coincide with $\so(1,n+1)$, and $\g$ must be
contained in $\simil(n)$. We show that the holonomy algebras $\g$
of types I and III do not annihilate any non-zero element in
$\nabla \R(\g)$, i.e. $\g$ must be of type II or IV. In this case
$(M,g)$ admits a parallel isotropic vector field, i.e. we reprove
the result of \cite{Sen08}. Next, we prove a reduction lemma that
allows to consider the following two cases: $\g=\h+p\wedge E$,
where $\h\subset\so(n)$ is an irreducible subalgebra, and
$\g=p\wedge E$.

We prove that the first case is impossible. For this we find the
form of all tensors in $\nabla \R(\g)$ annihilated by $\g$, it
turns out that this space is one-dimensional. Then we may find the
form of $\nabla R$. We calculate  $\nabla\Ric$, and show that the
Weyl conformal tensor $W$ is parallel ($\nabla W=0$). Then, using
the results of Derdzinski and Roter \cite{DR77,DR09} and of
\cite{GalConf}, we get a contradiction.

The second case corresponds to pp-waves \eqref{pp-wave}. The
condition $\nabla^2 R=0$ and simple computations allow us to find
the coordinate form of the metric.

The proof of this result from \cite{A-G} especially  shows the
power of the methods introduced in this paper, since lately there
appeared another more technical proof \cite{B-S-S} that uses
computations in local coordinates for metric \eqref{Walker}.

\section{Recurrent symmetric bilinear forms}\label{SecRecBF}

In \cite{Aminova} Aminova proved that if an indecomposable
Lorentzian manifold $(M,g)$ admits a parallel symmetric bilinear
form not proportional to the metric, then the manifold admits a
parallel isotropic vector field $p$, and the space of  parallel
symmetric bilinear forms is spanned over $\Real$ by the metric $g$
and by $\tau\otimes\tau$, where $\tau=g(p,\cdot)$ is the dual
1-form to $p$. We generalize this result to the case of recurrent
symmetric bilinear forms.

\begin{theorem} If an indecomposable Lorentzian
manifold $(M,g)$ admits a recurrent symmetric bilinear form not
proportional to the metric, then the manifold has holonomy
algebras contained in $\simil(n)$, in particular locally it is
given by the metric \eqref{Walker} and locally it admits recurrent
isotropic vector fields.

Let $g$ be given by \eqref{Walker} and suppose that it is
indecomposable. If $\partial_v^2H=\partial_i\partial_vH=0$, then
the coordinates can be chosen in such a way that $\partial_vH=0$,
in this case any recurrent symmetric bilinear form equals to
$f(\alpha g+\beta\tau\otimes\tau)$, where $\alpha,\beta\in\Real$,
$\tau=\d u=g(\partial_v,\cdot)$, and $f$ is a function; if
$\partial_v^2H\neq 0$ or $\partial_i\partial_vH\neq 0$, then any
recurrent symmetric bilinear is proportional either
 to $g$, or to $\tau\otimes\tau$.
\end{theorem}

For the proof it is enough  to find for the holonomy algebra
$\g\subset\simil(n)$ at a point $m\in M$ all invariant
1-dimensional subspaces in $\odot^2(\Real^{1,n+1})^*$ preserved by
$\g$. For algebras of type I and III these subspaces are $\Real
g_m$ and $\Real \tau_m\otimes \tau_m$. For algebras of type II and
IV these subspaces are 1-dimensional subspaces in $\Real g_m\oplus
\Real\tau_m\otimes \tau_m$. The condition
$\partial_v^2H=\partial_i\partial_vH=0$ holds only for the
holonomy algebras of type II and IV.

The above theorem can be used for studying Lorentzian manifolds
with recurrent Ricci tensors. We see that one deals with a Walker
metric, and the equations will be very similar to the Einstein
equation on the Walker metric, see \cite{GL10}.

\section*{Acknowledgements} I am thankful to D.V. Alekseevsky for
useful discussions and to Andrzej Derdzinski for helpful
correspondence.

\end{document}